\documentclass[reqno]{amsart}
\usepackage{verbatim}
\usepackage{amssymb}
\usepackage{enumerate}
\usepackage[active]{srcltx}

\usepackage[utf8x]{inputenc}
\usepackage[all]{xy}

\newtheorem{theorem}{Theorem}[section]

\newtheorem{lemma}[theorem]{Lemma}
\newtheorem{corollary}[theorem]{Corollary}

\newtheorem{shtheorem}{Shelah's Revised CGH Theorem\makebox[-2mm]{}}

   \usepackage{geometry}
\theoremstyle{definition}

\newtheorem{definition}[theorem]{Definition}

\theoremstyle{remark}

\newcommand{\mc}[1]{\mathcal{#1}}

\newcommand{\mb}[1]{\mathbf{#1}}

\newcommand{\setm}{\setminus}
\newcommand{\empt}{\emptyset}
\newcommand{\subs}{\subset}
\newcommand{\dom}{\operatorname{dom}}
\newcommand{\ran}{\operatorname{ran}}
\def\<{\left\langle}
\def\>{\right\rangle}
\def\br#1;#2;{\bigl[ {#1} \bigr]^ {#2} }

\newcommand{\kstr}{{\chi}}
\newcommand{\GG}{{\mc G}}

\newcommand{\bcal}{{\mc B}}
\newcommand{\pcal}{{\mc P}}
\newcommand{\bo}{{\beth_{\omega}}}
\newcommand{\acal}{{\mc A}}

\def\cf{\operatorname{cf}}
\newcommand{\MMM}{\mathbf M}

\newcommand{\MM}[3]{\MMM(#1,#2,#3)}

\newcommand{\CF}[4]{\ccf\big(#1,#2,#3\big)\le #4}

\newcommand{\ccf}{\operatorname{\mbox{${\chi}$}_{\rm CF}}}

\author[L. Soukup]{Lajos Soukup}
\thanks
  {Research  
supported by  the Hungarian National Research Grant OTKA  K 83726.} 
\address
      { Alfr{\'e}d R{\'e}nyi Institute of Mathematics, 
        Hungarian Academy of Sciences}
\email{soukup@renyi.hu}
\urladdr{http://www.renyi.hu/$\sim$soukup}
\keywords{revised GCH, almost disjoint, essentially disjoint, sparse, 
conflict-free coloring, conflict-free chromatic number, 
singular cardinal compactness}
\title[Essentially disjoint and  conflict free]
   {Essentially disjoint families, conflict free colorings and 
Shelah's Revised GCH}
\date{\today}
\begin{document}
\begin{abstract}Using Shelah's revised GCH theorem we prove that 
if ${\mu}<\beth_{\omega}\le \lambda$ are cardinals, then 
every ${\mu}$-almost disjoint family 
$\mc A\subs \br \lambda;\beth_{\omega};$ is {\em essentially disjoint},
i.e. for each $A\in \mc A$ there is a set
$F(A)\in \br A;<|A|;$ such that the family $\{A\setm
F(A):A\in \mc A\}$ is disjoint.

We also show that 
if  ${\mu}\le {\kappa}\le \lambda$ are  cardinals, ${\kappa}\ge {\omega}$,
and 
\begin{itemize}
 \item every ${\mu}$-almost disjoint family $\mc A\subs \br \lambda;{\kappa};$
is essentially disjoint,
\end{itemize}
then  
\begin{itemize}
 \item every  ${\mu}$-almost disjoint family 
$\mc B\subs \br \lambda;\ge {\kappa};$ has {\em a conflict-free coloring
with ${\kappa}$ colors}, i.e. there is a coloring $f:{\lambda}\to {\kappa}$
such that for all $B\in \mc B$ there is a color ${\xi}<{\kappa}$
such that $|\{{\beta}\in B: f({\beta})={\xi}\}|=1$.
\end{itemize}
Putting together these results we obtain that 
 if  ${\mu}<\beth_{\omega}\le \lambda$, then 
 every ${\mu}$-almost disjoint family 
$\mc B\subs \br \lambda;\ge \beth_{\omega};$ 
has a conflict-free coloring
with $\beth_{\omega}$ colors.

To yield  the above mentioned results we also need to prove 
a certain compactness theorem concerning singular cardinals.
\end{abstract}
\maketitle
\section{Introduction}

The investigation and comparison of properties of almost families of sets has 
a long history,
see e.g. the ancient paper of Miller, \cite{MIL}, the classical works
of Erdős, Galvin, Hajnal and Rothchild, \cite{EGH, EH1, EH2, EH3}, or 
the contributions of the next generation,  \cite{KO, KO2, KO3}
and \cite{Szeptycki}.

A set system $\mc A$ is {\em ${\mu}$-almost disjoint}
iff $|A\cap A'|<{\mu}$ for distinct  $A,A'\in \mc A$.

We say that $\mc A\subs \br X;{\kappa};$ is  {\em essentially
disjoint} ({\bf ED}, in short) iff for each $A\in \mc A$ there is a set
$F(A)\in \br A;<{\kappa};$ such that the family $\{A\setm
F(A):A\in \mc A\}$ is disjoint.

$\MM {\lambda}{\kappa}{\mu}\to
{\bf ED}$ denotes the statement that every ${\mu}$-almost disjoint
family $\mc A\subs \br {\lambda};{\kappa};$ is
  {\bf ED}.

In \cite[Theorems  2 and 5]{KO2} Komjáth  proved the following  results:
\begin{enumerate}[(A)]
 \item $\MM {\lambda}{\omega}r\to
{\bf ED}$ for $r<{\omega}$;
 \item $\MM {\lambda}{\omega_2}{\omega}\to
{\bf ED}$  under GCH;
 \item $\MM {\lambda}{{\omega}_1}{\omega}\to {\bf ED}
$  if GCH holds and ${\lambda}\le \aleph_{\omega}$;
 \item $\MM {\lambda}{{\omega}_1}{\omega}\to {\bf ED}
$  if V=L.
\end{enumerate}

Using Shelah' Revised GCH theorem,
in Section \ref{RGCH}
we will prove 
\begin{enumerate}[(A)]
\addtocounter{enumi}{4} 
 \item $\MM {\lambda}{{\beth_{\omega}}}{\mu}\to {\bf ED}
$ for each    ${\mu}<{\beth_{\omega}}$. 
\end{enumerate}

If  $f$ is a function  and  $A$ is a set, 
 we let
\begin{equation}\notag
I_f(A)=\{{\xi} \in \ran(f): |A \cap f^{-1}\{\xi\}|=1\}. 
\end{equation}
A function $f$ is a {\em conflict free coloring}
of a set system $\mc A$ iff $\dom(f)=\bigcup \mc A$ and    
$I_f(A)\ne \empt$ for all $A\in \mc A$.
\begin{equation}\notag
\CF {\lambda}{\kappa}{\mu}{\rho} 
\end{equation}
denotes the statement that every ${\mu}$-almost disjoint
family $\mc A\subs \br {\lambda};{\kappa};$ has a conflict-free coloring with 
${\rho}$ colors.

In \cite{HJSoSz} we proved that 
\begin{enumerate}[(a)]
 \item $\CF {\lambda}{\kappa}r{\omega}$ for $r<{\omega}\le {\kappa}\le {\lambda}$;
 \item $\CF {\lambda}{\kappa}{\omega}{{\omega}_2}$ 
for ${\omega_2}\le {\kappa}\le {\lambda}$ under GCH;
 \item $\CF {\lambda}{\kappa}{\omega}{{\omega}_1}$ for 
${\omega_1}\le {\kappa}\le {\lambda}\le \aleph_\omega$  under GCH;
 \item $\CF {\lambda}{\kappa}{\omega}{{\omega}_1}$ for 
${\omega_1}\le {\kappa}\le {\lambda}$ if V=L.
\end{enumerate}

In \cite[Corollary 2]{KO3}   Komjáth improved (a) by showing
\begin{enumerate}[(a${}^*$)]
 \item $\CF {\lambda}{\ge {\omega}}r{{\omega}}$  for $r<{\omega}\le {\kappa}$,
\end{enumerate}
i.e
 every $r$-almost disjoint family 
$\mc A\subs \br {\lambda};\ge {\omega};$, where $r\in {\omega}$,   
has a conflict free
coloring with ${\omega}$ colors. 

In this paper 
we will show 
\begin{enumerate}[(a${}^*$)]
\addtocounter{enumi}{1} 
\item 
$\CF {\lambda}{\ge{\omega}_2}{\omega}{{\omega}_2}$ under GCH;
\item $\CF {\lambda}{\ge{\omega}_1}{\omega}{{\omega}_1}$ if GCH holds
and ${\lambda}\le \aleph_{\omega}$;
\item $\CF {\lambda}{\ge{\omega}_1}{\omega}{{\omega}_1}$ if V=L;
 \item $\CF {\lambda}{\ge{\beth_{\omega}}}{\nu}{\beth_{\omega}}$ 
for each   ${\nu}<{\beth_{\omega}}$. 
\end{enumerate}

We do not need to present four proofs, because 
 the following implication holds between essential disjointness
and conflict free colorings:
\begin{equation}
\MM {\lambda}{\kappa}{\mu}\to
{\bf ED}\text { implies } \CF {\lambda}{\ge \kappa}{\mu}{{\kappa}}, 
\end{equation}
see Corollary \ref{cor:main}.

To prove our results, in section \ref{sc:pmscc} we formulate a 
theorem which is a poor relative of Shelah's 
Singular Cardinal Compactness Theorem 
(\cite[Theorem 2.1]{Sh-SCC}).

Our notation is standard, see e.g.  \cite{KU}.

If $X$ is a set and ${\kappa}$ is a cardinal, then  write
$\br X;\ge {\kappa};=\{A\subs X: |A|\ge {\kappa}\}$.

If $\<X_{\alpha}:{\alpha}<\sigma\>$ is a sequence of sets, we will often write
$X_{< {\alpha}}=\bigcup_{{\beta}< {\alpha}}X_{\alpha}$, $X_{\le {\alpha}}=\bigcup_{{\beta}\le {\alpha}}X_{\alpha}$, etc.

\section{Poor man's singular cardinal compactness}\label{sc:pmscc}

Both the formulation and the proof of the  following statement
use  the ideas of Shelah's Singular Cardinal Compactness theorem 
(see \cite[Theorem 2.1]{Sh-SCC}).

A set system $\GG$  is {\em $\kstr$-chain closed}
iff $\bigcup_{{\alpha}<\kstr}G_{\alpha}\in \GG$
for  any $\subseteq$-increasing sequence $\<G_{\alpha}:{\alpha}<\kstr\>\subs \GG$.

\begin{theorem}
\label{lm:compact}
Assume that ${\lambda}$ is a singular cardinal
and $\GG\subs \br {\lambda};<{\lambda};$.
If for some cardinals $\kstr\le {\mu}<{\lambda}$,  
\begin{enumerate}[{$(\circ)$}]
\item \quad the set $\GG \cap \br {\lambda};{\nu};$ is {$\kstr$-chain closed}
   and  cofinal in $\br {\lambda};{\nu};$
for each ${\mu}\le {\nu}<{\lambda}$,
\end{enumerate}
then  
there is a continuous, increasing sequence
$\<G_{\xi}:{\xi}<\cf({\lambda})\>\subs \GG$ such that  
$\bigcup_{{\xi}<\cf({\lambda})}G_{\xi}={\lambda}$.    
\end{theorem}

\begin{proof}[Proof of Theorem \ref{lm:compact}]
For each $Y\in  {\lambda}$ with ${\mu}\le |Y|<{\lambda}$
by (b) we can pick $G(Y)\in \GG\cap \br {\lambda};|Y|;$ with $G(Y)\supset Y$.

Let $\<{\lambda}_{\zeta}:{\zeta}<\cf({\lambda})\>$ be a strictly increasing, 
continuous and cofinal sequence of cardinals in ${\lambda}$ with  
${\lambda}_0>\mu,\kstr, \cf({\lambda})$.

By transfinite induction on $n<\kstr$
we will define sets $$\<Y_{{\zeta},n} :{\zeta}<\cf({\lambda}),
n<\kstr \>\subs \br {\lambda};<{\lambda};$$
such that 
\begin{enumerate}[(A)]
\item $|Y_{{\zeta},n}|={\lambda}_{\zeta}$,
 \item  
the sequences $\<Y_{{\zeta},n} :{\zeta}<\cf({\lambda})\>$ are
increasing and continuous, 
\item $G\big(\bigcup_{m<n}Y_{{\zeta},m}\big)\subs Y_{{\zeta},n}$,
\end{enumerate}
as follows (see Figure 1).

Let $Y_{{\zeta},0}={\lambda}_{\zeta}$ for ${\zeta}<\cf({\lambda})$.

Assume that the family $$\<Y_{{\zeta},m} :{\zeta}<\cf({\lambda}),
m<n \>\subs \br {\lambda};<{\lambda};$$ is defined.

Let 
\begin{equation}\notag
B_{{\zeta},n}=G\big(\bigcup_{m<n}Y_{{\zeta},m}\big),
\end{equation}
and write
\begin{equation}\notag
B_{{\zeta},n}=\{b_{{\zeta},n}(i):i<{\lambda}_{\zeta}\} 
\end{equation}
for ${\zeta}<\cf({\lambda})$, and for ${\xi}<\cf({\lambda})$ let
\begin{equation}\notag
Y_{{\xi},n}=\{b_{{\zeta},n}(i):{\zeta}<\cf({\lambda}), i<{\lambda}_{\xi}\}. 
\end{equation}

\begin{figure}
\begin{displaymath}
\xymatrix{
{\kstr}&{G_{0}}\ar@{}[r]|*{\subset} &{G_{\zeta}} & & &
\\
&& & & & &\\&*[r]{Y_{0,1}}\ar@{}[r]|*{\subs}\ar@{}[u]|{\quad \bigcup}
&*[r]{Y_{{\zeta},1}}\ar@{.>}^{continuous}[rrrr]\ar@{}[u]|{\quad \bigcup}& & & &
\\
&*[r]{G(Y_{0,0})\in \GG}\ar@{}[u]|{\quad \bigcup}&*[r]{G(Y_{{\zeta},0})\in \GG}\ar@{}[u]|{\quad \bigcup}& & & &\\  
&*[r]{Y_{0,0}}\ar@{}[u]|{\quad \bigcup}
\ar@{}[r]|*{\dots}&*[r]{Y_{{\zeta},0}}\ar@{}[u]|{\quad
  \bigcup}\ar@{.>}^{continuous}[rrrr]&&
&&\\  
\ar[uuuuu]
\ar[rrrrrr]&&&&&&&\cf({\lambda})
}
\end{displaymath}
\caption{}
\end{figure}

Since the sequence 
$\<{\lambda}_{\zeta}:{\zeta}<\cf({\lambda})\>$
was continuous, the sequence 
$\<Y_{{\zeta},n} :{\zeta}<\cf({\lambda})\>$ is also continuous, so (B) holds.

(A) and (C) are clear from the construction
because $\cf({\lambda})<{\lambda}_0$.

Let 
\begin{equation}\notag
 G_{{\xi}}=\bigcup_{n<\kstr }Y_{{\xi},n}.
\end{equation}
for ${\xi}<\cf({\lambda})$.

To check $G_{\xi}\in \GG$ it is enough to observe that 
\begin{equation}\notag
 G_{{\xi}}=\bigcup_{n<\kstr }G(Y_{{\xi},n}).
\end{equation}
because 
the sequence 
$\<G(Y_{{\xi},n}):n<\kstr\>\subs\GG\cap \br {\lambda};{\lambda}_{\xi};$ is $\subseteq$-increasing and  
$\GG\cap \br {\lambda};{\lambda}_{\xi};$  is $\kstr$-chain closed.

By (B) we have   $Y_{{\zeta},n}\subs Y_{{\xi},n}$
for ${\zeta}<{\xi}<\cf({\lambda})$, and so  $G_{\zeta}\subs G_{\xi}$.

Since ${\lambda}_{\zeta}\subs G_{\zeta}$, we also have 
$\bigcup_{{\xi}<\cf({\lambda})}G_{\xi}={\lambda}$.

Finally assume that ${\xi}<\cf({\lambda})$ is a limit ordinal.
Then
\begin{equation}\notag
 G_{\xi}=\bigcup_{n<\kstr}Y_{{\xi},n}=\bigcup_{n<\kstr}\Big(
\bigcup_{{\zeta}<{\xi}}Y_{{\zeta},n}\Big)=
\bigcup_{{\zeta}<{\xi}}\Big(
\bigcup_{n<\kstr}Y_{{\zeta},n}\Big)=
\bigcup_{{\zeta}<{\xi}}G_{\zeta},
\end{equation}
so the sequence  $\<G_{\xi}:{\xi}<\cf({\lambda})\>$
is continuous, which was to be proved.
\end{proof}

\section{Families of sets of size  $\beth_{\omega}$}\label{RGCH}

\begin{definition}
If $\nu\le \rho$ are cardinals, then write  
\begin{equation}\notag
{\rho}^{[{\nu}]}={\rho} 
\end{equation}
iff there is a family
$\bcal\subs \br {\rho};\le {\nu};$ of size ${\rho}$
such that for all $u\in \br {\rho};{\nu};$
there is $\pcal\in \br \bcal;<{\nu};$ such that 
$u\subs \cup\pcal$.   
\end{definition}

\begin{shtheorem}[{\cite[Theorem 0.1]{Sh-GCH}}]
\label{tm:sh}
\mbox{} \\If ${\rho}\ge \bo$, then ${\rho}^{[{\nu}]}={\rho}$
for each large enough regular cardinal ${\nu}<\bo$.  
\end{shtheorem}

\begin{lemma}\label{lm:bo_closed}Fix ${\mu}<\bo$. Then for each 
${\rho}\ge \bo$
there is a  regular ${\nu}({\rho})<\bo$ 
such that if 
  $\acal\subs \br \rho;{\nu}(\rho);$ is
${\mu}$-almost disjoint, then $|\acal|\le \rho$.
\end{lemma}

\begin{proof}[Proof of lemma \ref{lm:bo_closed}]
Let ${\mu}<{\nu}<\bo$ be regular such that
${\rho}^{[{\nu}]}={\rho}$ witnessed by a family 
  $\bcal\subs \br \rho;{\nu};$.   
We show that  $\nu(\rho)=\nu$ works.

If $A\in \acal$, 
then there is $B\in \bcal$ such that $A\cap B={\nu}$.
If $A_0,A_1\in \acal$ are distinct, and 
$|A_0\cap B|=|A_1\cap B|={\nu}$ then 
$A_0\cap B\ne A_1\cap B$ because 
$\acal$ is ${\mu}$-almost disjoint.
So
\begin{displaymath}
 |\{A\in \acal:|A\cap B|={\nu}\}|\le 2^{|B|}\le 2^{\nu}<\bo. 
\end{displaymath}
Thus $  |\acal|\le \bo\cdot|\bcal|=\rho$.

\end{proof}

\begin{lemma}\label{lm:loc_small}
If ${\lambda}> {\kappa}\ge \bo >{\mu}$, and 
$\{A_{\alpha}:{\alpha}<\tau\}\subs \br {\lambda};{\kappa};$
is a ${\mu}$-almost disjoint family, then $\tau\le {\lambda}$, and  
 there is an increasing, continuous sequence 
$\<G_{\zeta}:{\zeta}<\cf({\lambda})\>\subs \br {\lambda};<{\lambda};$
 such that 
\begin{equation}\label{eq:loc_small}
\text{
$\forall {\zeta}<\cf({\lambda})$\
$\forall {\alpha}\in G_{\zeta+1}\setm G_{\zeta}$\ $(\ |A_{\alpha}\cap G_{\zeta}|<\bo$
and $A_{\alpha}\subs G_{{\zeta}+1}\ )$.}
\end{equation}
\end{lemma}

\smallskip
\noindent {\em Remark.} In the published version of this paper we considered only 
the special case when ${\kappa}=\bo$.

\begin{proof}
Applying lemma \ref{lm:bo_closed} for $\rho={\lambda}$ we have 
$|\tau|\le {\lambda}$.
So we can assume that $\tau={\lambda}$. 
 
We should distinguish two cases. 

\medskip\noindent{\bf Case 1.}
{\em ${\lambda}$   is a regular cardinal}

\newcommand{\card}{\mb{Card}} 
Denoting  by $\card$ the class of cardinals,
pick  a cardinal ${\mu}<{\nu}<\bo$ such that
\begin{displaymath}
  B=\{{\rho}\in \card\cap {\lambda}:{\nu}({\rho})={\nu}\}
\end{displaymath}
is cofinal in $\card\cap {\lambda}$.
(If ${\lambda}={\sigma}^+$   is a successor cardinal,
then  $B=\{\sigma\}$ and $\nu=\nu(\sigma)$  work.)

Let 
\begin{equation}\notag
B^*=\{{\zeta}<{\lambda}:|{\zeta}|\in B\},
\end{equation}
and for $\zeta\in B^*$ let
\begin{equation}\notag
 f({\zeta})=A_{<{\zeta}}\cup \{{\alpha}: |A_{\alpha}\cap {\zeta}|\ge {\nu}\}.
\end{equation}
Since the family 
$\{A_{\alpha}\cap {\zeta}:  |A_{\alpha}\cap {\zeta}|\ge {\nu}\}
\subs \br {\zeta};\ge {\nu};$ is ${\mu}$-almost disjoint, 
by Lemma \ref{lm:bo_closed} we have $|f({\zeta})|\le |\zeta|$.
Let 
\begin{equation}\notag
 D=\{{\delta}<{\lambda}:\cf({\delta})=\nu^+,\ \sup(B^*\cap {\delta})={\delta}\ \land\  
(\forall {\zeta}\in B^*\cap {\delta}) 
\ f({\zeta})\subs {\delta}\}.
\end{equation}
Since $B^*$ is cofinal in ${\lambda}$, ${\nu}^+\le \bo<{\lambda}$
and $f:B^*\to \br {\lambda};<{\lambda};$, the set
$D$ is also cofinal in ${\lambda}$, and 
\begin{equation}\notag
   \forall {\delta}\in D(\ \forall {\alpha}<{\delta}
\ A_{\alpha}\subs {\delta} \land \forall {\alpha}\in ({\lambda}\setm {\delta})
\ |A_{\alpha}\cap {\delta}|<\nu).
\end{equation}
Indeed, if $|A_{\alpha}\cap {\delta}|\ge{\nu}$, then there is 
${\zeta}\in B^*\cap {\delta}$ with $|A_{\alpha}\cap \delta|\ge \nu$ 
for $\cf({\delta})={\nu}^+$,
and  then ${\alpha}\in f({\zeta})\subs \delta$ by the definition of $f$.

Thus we also have  
\begin{equation}\label{eq:d'}
   \forall {\gamma}\in D'(\ \forall {\alpha}<{\gamma}
\ A_{\alpha}\subs {\gamma} \land \forall {\alpha}\in ({\lambda}\setm {\gamma})
\ |A_{\alpha}\cap {\gamma}|<\bo). 
\end{equation}
Indeed,  
if $|A_{\alpha}\cap {\gamma}|=\bo$, then there is 
$\delta\in D\cap (\gamma+1)$ with $|A_{\alpha}\cap \delta|\ge \nu$,
and so $A_{\alpha}\subs \delta$.

Let $\{\gamma_{\zeta}:{\zeta}<{\lambda}\}$ be the increasing enumeration
of the club set $D'$.
Then the choice $G_{\zeta}=\gamma_{\zeta}$ works
by (\ref{eq:d'}).

\medskip\noindent{\bf Case 2.}
{\em ${\lambda}>\cf({\lambda})$   is a singular  cardinal.}

Let 
\begin{equation}\notag
 \GG=\{G\in \br {\lambda};<{\lambda};:|G|\ge{\kappa}\land (\forall {\alpha}\in G)\ 
A_{\alpha}\subs G \land  (\forall {\alpha}\in {\lambda}\setm G)\   
|A_{\alpha}\cap G|<\nu(|G|\}.
\end{equation}

If ${\kappa}\le \rho<{\lambda}$, then 
${\nu}(\rho)>{\omega}$ is a regular cardinal, so the family 
$\GG\cap \br {\lambda};\rho;$ is ${\omega}$-chain closed.

Next we show that the set $\GG\cap \br {\lambda};\rho;$ is cofinal in 
$\br {\lambda};\rho;$
for all ${\kappa}\le \rho <{\lambda}$.

Indeed, let $Y\in \br {\lambda};\rho;$.
Define an increasing sequence $\<Y_n:n<{\omega} \>\subs \br {\lambda};\rho;$ 
as follows:
\begin{enumerate}[(i)]
 \item $Y_0=Y$,
\item for $1\le n<{\omega}$ let 
\begin{equation}\notag
Y_n=Y_{<n}\cup
\bigcup \{A_{\alpha}:{\alpha}\in Y_{n-1}\}\cup
\{{\alpha}:|Y_{<n}\cap A_{\alpha}|\ge \nu(\rho)\}.
\end{equation}
\end{enumerate}
We are to show that $G=\bigcup_{n<{\omega}}Y_n\in \mc G\cap \br {\lambda};\rho;$.

By induction on $n$ we obtain $|Y_n|=\rho$ because
$|\{{\alpha}:|Y_{<n}\cap A_{\alpha}|\ge \nu(\rho)\}|\le \rho$ 
by lemma \ref{lm:bo_closed}.
So $|G|=\rho$.

If ${\alpha}\in G$, then ${\alpha}\in Y_n$ for some $n$, and so 
$A_{\alpha}\subs Y_{n+1}\subs G$.

If $|A_{\alpha}\cap G|\ge {\nu}(\rho)$, then $|A_{\alpha}\cap Y_n|\ge {\nu}(\rho)$
for some $n<{\omega}$ because  $\cf({\nu}(\rho))={\nu}(\rho)>{\omega}$.
So ${\alpha}\in Y_{n+1}\subs G$. 

Thus we proved $G\in \mc G$ and so $\GG\cap \br {\lambda};\rho;$ is really 
cofinal in 
$\br {\lambda};\rho;$.

So we can apply Theorem \ref{lm:compact} to obtain the required
sequence  $\<G_{\zeta}:{\zeta}<\cf({\lambda})\>$.
\end{proof}

\begin{theorem}\label{tm:bo}
$\MM {\lambda}{{\beth_{\omega}}}{\mu}\to {\bf ED}
$ whenever    ${\nu}<{\beth_{\omega}}\le {\lambda}$. 
\end{theorem}

\begin{proof}We prove the theorem 
by induction on ${\lambda}$.
If ${\lambda}=\bo$, then 
applying lemma \ref{lm:bo_closed} for $\rho=\bo$ 
we have $|\acal|\le \bo$, so 
the statement is clear.

Assume now that ${\lambda}>\bo$ and we proved the theorem 
for ${\lambda}'<{\lambda}$.

Let $\mc A\subs \br {\lambda};\bo;$ be ${\mu}$-almost disjoint.
By lemma \ref{lm:bo_closed}, $|\mc A|\le {\lambda}$, so we can write
$\mc A=\{A_{\alpha}:{\alpha}<{\lambda}\}$.

By lemma  \ref{lm:loc_small} there is an increasing, continuous sequence 
$\{G_{\zeta}:{\zeta}<\cf({\lambda})\}\subs \br {\lambda};<{\lambda};$
which satisfies the  (\ref{eq:loc_small}).

For ${\zeta}<\cf({\lambda})$ write
\begin{equation}\notag
 X_{\zeta}=G_{{\zeta}+1}\setm G_{\zeta}\text{ and }
\mc A_{\zeta}=\{A_{\alpha}\cap X_{\zeta}:{\alpha}\in X_{\zeta}\}.
\end{equation}
Then $\mc A_{\zeta}\subs \br X_{\zeta};\bo;$ is ${\mu}$-almost disjoint,
so by the inductive assumption 
there are set $F_{\zeta}(A_{\alpha})\in \br A_{\alpha};<\bo;$
for ${\alpha}\in X_{\zeta}$ such that the family
\begin{equation}\notag
 \{(A_{\alpha}\cap X_{\zeta})\setm F_{\zeta}(A_{\alpha}):{\alpha}\in X_{\zeta}\}
\end{equation}
is disjoint.
Let
\begin{equation}\notag
 F(A_{\alpha})=F_{\zeta}(A_{\alpha})\cup (A_{\alpha}\cap X_{<{\zeta}})
\end{equation}
for ${\alpha}\in X_{\zeta}$.
Then  $|F(A_{\alpha})|<\bo$, and   the sets 
\begin{equation}\notag
 \{A_{\alpha}\setm F(A_{\alpha}):{\alpha}<{\lambda}\}
\end{equation}
are disjoint.

\end{proof}

\section{Essential disjointness and conflict free colorings}

We say that a set system $\mc A$ is  
{\em ${\kappa}$-hereditary essentially disjoint, }
({\bf ${\kappa}$-hED}, in short)
iff for each function $H$ with $\dom(H)=\mc A$ and $H(A)\in \br A;{\kappa};$
for $A\in \mc A$,
the family $\{H(A):A\in \mc A\}$ is {\bf ED}.

A ${\kappa}$-{\bf hED} family is  clearly 
${\kappa}$-almost disjoint.

\begin{theorem}\label{tm:main}
Let ${\lambda}\ge {\kappa}\ge {\omega}$ be  cardinals and $\mc A\subs \br {\lambda};\ge {\kappa};$
be a family of sets.
If 
\begin{enumerate}[(a)]
 \item[$\bullet$] $\mc A$ is ${\kappa}$-hereditarily essentially disjoint,
\end{enumerate}
then 
\begin{enumerate}[(a)]
 \item[$\circ$] there is a coloring $c:{\lambda}\to {\kappa}$
such that $|I_c(A)|={\kappa}$ for all $A\in \mc A$.
\end{enumerate}

\end{theorem}

\begin{corollary}\label{cor:main}
Let ${\mu}\le {\kappa}\le {\lambda}$ be cardinals, ${\kappa}\ge {\omega}$.
If 
\begin{enumerate}[$(*_{\lambda})$]
 \item every ${\mu}$-almost disjoint family $\mc A\subs \br \lambda;{\kappa};$
is {\bf ED},
\end{enumerate}
then  
\begin{enumerate}[$(\star_{\lambda})$]
 \item for each  ${\mu}$-almost disjoint family 
$\mc A\subs \br \lambda;\ge {\kappa};$
there is a coloring $f:\lambda\to {\kappa}$
such that $|I_f(A)|={\kappa}$ for all $A\in \mc A$,
and so $\CF {\lambda}{\ge\kappa}{\mu}{\kappa}$.
\end{enumerate}
\end{corollary}

In the inductive proof of 
Theorem \ref{tm:main}  we will use the following 
observation on decomposability   of ${\kappa}$-{\bf hED} families.

\begin{theorem}\label{tm:decomposition}
If ${\lambda}>{\kappa}\ge {\omega}$ are cardinals, and 
$\mc A=\{A_{\alpha}:{\alpha}<{\lambda}\}\subs \br {\lambda};\ge {\kappa};$ is 
a ${\kappa}$-{\bf hED}  family, then 
there is a  partition $\{X_{\zeta}:{\zeta}<\cf({\lambda})\}
\subs \br {\lambda};<{\lambda};$
of ${\lambda}$
such that 
\begin{enumerate}[(M1)]
\item $A_{\alpha}\cap A_{{\alpha}'}\subs  X_{\le \zeta}$ 
for distinct ${\alpha},{\alpha}'\in X_{\le\zeta}$.
 \item $|A_{\alpha}\cap X_{<\zeta}|<{\kappa}$
for all ${\alpha}\in X_{\zeta}$,
\item 
$|A_{\alpha}\cap X_{\zeta}|\ge {\kappa}$ for all ${\alpha}\in  X_{\zeta}$,
\end{enumerate}
\end{theorem}

\begin{proof}[Proof of Theorem \ref{tm:decomposition}]
We should distinguish two cases.

\medskip
\noindent{\bf Case 1. ${\lambda}$ is regular.}

For each $A\in \mc A$, let $H(A)$ be the first ${\kappa}$ many elements of $A$.
Since the family $\{H(A):A\in \mc A\}$ is $\bf ED$, there is an 
 injective function $f$  on $\mc A$ with $f(A)\in H(A)$.

Fix a large enough regular cardinal $\theta$, and let 
$\<N_{\zeta}:1\le {\zeta}<{\lambda}\>$ 
be a strictly increasing continuous sequence of elementary 
submodels of $\<H_{\theta},\in \>$ such that 
\begin{enumerate}[(i)]
 \item ${\kappa}+{\zeta}\subs N_{\zeta}\cap {\lambda}\in {\lambda}$ and $|N_{\zeta}|={\kappa}+|{\zeta}|$, 
\item $\<A_{\alpha}:{\alpha}<{\lambda}\>, H, f\in N_1$.
\item $N_{\zeta}\in N_{{\zeta}+1}$.
\end{enumerate}
Write $N_0=\empt$.

For  ${\zeta}<cf({\lambda})$ let
\begin{displaymath}
X_{\zeta}=(N_{{\zeta}+1}\setm N_{\zeta})\cap {\lambda}.
\end{displaymath}

 (M1) is  clear because $N_{{\zeta}+1}$ is an elementary submodel,
$X_{\le {\zeta}}=N_{\zeta+1}\cap {\lambda}$, 
and ${\kappa}+1\subs N_{{\zeta}+1}$.

To check (M2), assume that ${\alpha}<{\lambda}$ with 
$|A_{\alpha}\cap X_{<{\zeta}}|={\kappa}$.
Since $H(A_{\alpha})$ is the first ${\kappa}$ many elements of $A_{\alpha}$,
and $X_{<{\zeta}}=N_{\zeta}\cap {\lambda}$
is an initial segment of ${\lambda}$, we have $H(A_{\alpha})\subs N_{\zeta}$.
So ${\eta}=f(A_{\alpha})\in H(A_{\alpha})\subs N_{\zeta}$.
But ${\alpha}$ is definable from $f$, $\<A_{\alpha}:{\alpha}<{\lambda}\>$ and ${\eta}$, because $f$ was injective.
So ${\alpha}\in N_{\zeta}\cap {\lambda}=X_{<{\zeta}}$.

To check (M3) assume that 
${\alpha}\in (\mc N_{\zeta+1}\setm N_{\zeta})\cap {\lambda}$.
Then $|A_{\alpha}\cap N_{{\zeta}+1}|\ge {\kappa}$ because 
$N_{{\zeta}+1}$ is an elementary submodel
and ${\kappa}+1\subs N_{{\zeta}+1}$.
On  the other hand, $|A\cap N_{{\zeta}}|<{\kappa}$ by (M2),
so $|A_{\alpha}\cap ( N_{\zeta+1}\setm N_{\zeta})|={\kappa}$.

\medskip
\noindent{\bf Case 2. ${\lambda}$ is singular.}

Let 
\begin{equation}\notag
 \kstr=\left\{
\begin{array}{ll}
{\omega}&\text{$cf({\kappa})\ne {\omega}$,}\bigskip\\
{\omega}_1&\text{$cf({\kappa})={\omega}$.}
\end{array}
\right .
\end{equation}
We say that $Y\subs {\lambda}$ is {\em good}
if (G1)-(G3) below hold:
\begin{enumerate}[(G1)]
\item if ${\alpha}\ne{\beta}\in Y$, then 
$A_{\alpha}\cap A_{\beta}\subs Y$,
\item  if ${\alpha}\in Y$ then 
$|A_{\alpha}\cap  Y|\ge {\kappa}$,
\item  if  
$|A_{\alpha}\cap  Y|\ge {\kappa}$  then ${\alpha}\in Y$.
\end{enumerate}

Let 
\begin{equation}\notag
 \GG=\{Y\in \br {\lambda};<{\lambda};:\text{$Y$ is good}\}.
\end{equation}

\begin{lemma}
$\GG\cap \br {\lambda};{\nu};$ is $\kstr$-chain closed for 
$\max(\kappa,\cf({\lambda}), \kstr)< {\nu}<{\lambda}$. 
\end{lemma}

\begin{proof}
Assume that  $\<Y_n:n<\kstr \>\subs \GG\cap \br {\lambda};{\nu};$ is increasing.

(G1) and  (G2) are clear.

To check (G3) assume that $|A_{\alpha}\cap \bigcup_{n<\kstr}Y_n|\ge{\kappa}$.
Since $\cf({\kappa})\ne \kstr$, there is $n<\kstr$ with 
$|A_{\alpha}\cap Y_n|={\kappa}$, and so ${\alpha}\in Y_{n}$.
\end{proof}

\begin{lemma} 
$\GG\cap \br {\lambda};{\nu};$ is cofinal in $\br {\lambda};{\nu};$ for 
$\max(\kappa,\cf({\lambda}), \kstr)< {\nu}<{\lambda}$.

\end{lemma}

\begin{proof}[Proof of the Lemma]
Let  $Y\in \br {\lambda};{\nu};$.

For each $A\in \mc A$ pick  $H(A)\in \br A;{\kappa};$.

Define an increasing sequence $\<Y_n:n<\kstr \>\subs \br {\lambda};{\nu};$ 
as follows:
\begin{enumerate}[(i)]
 \item $Y_0=Y$,
\item for $1\le n<\kstr$ let 
\begin{multline}\notag
Y_n=Y_{<n}\cup
\bigcup\{A_{\alpha}\cap A_{\beta}:{\alpha}\ne {\beta}\in Y_{<n}\}\cup
\\
\bigcup\{H(A_{\alpha}):{\alpha}\in Y_{<n}\}\cup
\{{\alpha}:|Y_{<n}\cap A_{\alpha}|\ge {\kappa}\}.
\end{multline}

\end{enumerate}
Let $G=\bigcup_{n<\kstr}Y_n$.
We are to show that $G\in \mc G\cap \br {\lambda};{\nu};$.

Since $|Y_n|\le |Y_{<n}|+{\kappa}$ by  inequality (\ref{lm:small_trace}), 
we have ${\nu}\le |Y|\le |G|\le |Y|+{\kappa}+\kstr=
{\nu}$.

Properties 
(G1) and (G2) are straightforward from the construction.

To check (G3) assume that $|A_{\alpha}\cap G|\ge{\kappa}$.
Since $\cf({\kappa})\ne \kstr$, there is $n<\kstr$ with 
$|A_{\alpha}\cap Y_n|={\kappa}$, and so ${\alpha}\in Y_{n+1}\subs G$.
\end{proof}

By Theorem \ref{lm:compact} there is a continuous, increasing sequence
$\<G_{\xi}:{\xi}<\cf({\lambda})\>\subs \GG$ such that  
$\bigcup_{{\xi}<\cf({\lambda})}G_{\xi}={\lambda}$.

Now for ${\zeta}<\cf({\lambda})$ let
\begin{equation}\notag
 X_{\zeta}=G_{\zeta+1}\setm G_{\zeta}. 
\end{equation}

Clearly  $|X_{\zeta}|\le |G_{\zeta}|<{\lambda}$.
The sequence  
$\<G_{{\xi}}:{{\xi}}<\cf({\lambda})\>$ was continuous, so
the family  $\{X_{\zeta}:{\zeta}<\cf({\lambda})\}$ is a partition of
 ${\lambda}$. 

Since $G_{\zeta}\in \GG$ we have 
\begin{equation}\label{eq:ay}
  |A_{\alpha}\cap G_{\zeta}|\ge {\kappa} \text{ iff } {\alpha}\in G_{\zeta}. 
\end{equation}

(M1) is clear because 
$X_{\le {\zeta}}=G_{{\zeta}+1}\in \GG$ and (G1) holds for $G_{{\zeta}+1}$ .

To check (M2) assume that 
${\alpha}\in X_{\zeta}$,
and so ${\alpha}\notin G_{\zeta}$. Thus 
$|A_{\alpha}\cap G_{\zeta}|<{\kappa}$ by (\ref{eq:ay}).
But $G_{\zeta}=X_{<{\zeta}}$, so (M2) holds.

To check (M3) assume that ${\alpha}\in X_{\zeta}$.
Then ${\alpha}\in G_{{\zeta}+1}\setm G_{\zeta}$ and so applying (\ref{eq:ay})
twice we obtain
\begin{equation}\notag
   |A_{\alpha}\cap G_{{\zeta}+1}|\ge {\kappa} \land 
   |A_{\alpha}\cap G_{{\zeta}}|<{\kappa}. 
\end{equation}
So $|A_{\alpha}\cap X_{\zeta}|\ge {\kappa}$.

So we proved Theorem \ref{tm:decomposition}.
\end{proof}

\begin{proof}[Proof of Theorem \ref{tm:main}]
 Let ${\mc I}$ be a principal ideal on ${\kappa}$
with $\br {\kappa};<{\kappa};\subs \mc I$.
We prove the following stronger statement:

\begin{itemize}
  \item[($\bullet_{\lambda}$)]  
If a family $\mc A\subs \br {\lambda};\ge {\kappa};$
is { ${\kappa}$-{\bf hED}},
and   $F:{\lambda}\to \mc I$, 
then  there is a function $c:{\lambda}\to {\kappa}$  such that 
\begin{enumerate}[({$\mc I$}1)]
\item $c({\xi})\notin F({\xi})$ for ${\xi}\in {\lambda}$,
\item  $|I_c(A)|={\kappa}$ for all $A\in {\mc A}$. 
\end{enumerate}
\end{itemize}

\medskip

If  $Y\subs {\lambda}$, then  consider the family 
\begin{equation}\notag
\mc A\lceil Y=\{A\in \mc A: |A\cap Y |\ge {\kappa}\}, 
\end{equation}
and for each $A\in \mc A\lceil Y$ pick $H(A)\in \br A\cap Y;{\kappa};$.
Since $\mc A$ is ${\kappa}$-{\bf hED}, the family 
$\{H(A):A\in \mc A\lceil Y\}\subs \br Y;{\kappa};$ is 
 essentially disjoint, so 
\begin{equation}\label{lm:small_trace}
|\mc A\lceil Y|\le |Y|. 
\end{equation}

Especially, $|\mc A|\le {\lambda}$, and so we can  
write $\mc A=\{A_{\alpha}:{\alpha}<{\lambda}\}$.

We prove ($\bullet_{\lambda}$) by induction on ${\lambda}$.

\medskip

If ${\lambda}= {\kappa}$, then  let $c:{\lambda}\to {\kappa}$ be 
any injective function such that $c({\xi})\notin F({\xi})$ for 
${\xi}\in {\lambda}$.
 Then $I_{c}(A)=c''A\in \br {\kappa};{\kappa};$ for all 
$A\in \mc A$. 

\medskip

Assume now that 
${\lambda}>{\kappa}$ and 
$(\bullet_{\lambda'})$ holds for ${\lambda}'<{\lambda}$.

Apply Theorem \ref{tm:decomposition} for the family 
$\{A_{\alpha}:{\alpha}<{\lambda}\}$
to find a partition $\{X_{\zeta}:{\zeta}<\cf({\lambda})\}
\subs \br {\lambda};<{\lambda};$
of ${\lambda}$ satisfying (M1)--(M3).

Let $\mc A_{\zeta}=\{A_{\alpha}:{\alpha}\in X_{\zeta}\}\in \br \mc A;<{\lambda};$ for 
${\zeta}<\cf({\lambda})$. 

By recursion, for each  ${\zeta}<\cf({\lambda})$ we define
a function \begin{displaymath}
 c_{\zeta}: X_{\zeta}\to {\kappa}
\end{displaymath}
and a function
$G_{\zeta}:\mc A_{\zeta}\to {\mc I}\cap \br {\kappa};{\kappa};$
such that 
\begin{equation}\label{eq:GsubsI}
 G_{\zeta}(A_{\alpha})\subs I_{c_{\le{\zeta}}}(A_{\alpha}) 
\text{ for all ${\alpha}\in X_{\zeta}$}
\end{equation}
as follows. 

Assume that we have defined $\<c_{\eta}:{\eta}<{\zeta}\>$
and $\<G_{\eta}:{\eta}<{\zeta}\>$.

If $x\in X_{\zeta}$ then there is at most one 
${\alpha}\in X_{<\zeta}$ with $x\in A_{\alpha}$ by (M1).
Define $F_{\zeta}:X_{\zeta}\to {\mc I}$ as follows:
\begin{displaymath}
 F_{\zeta}(x)=\left\{
 \begin{array}{ll}
F(x)\cup G_{\eta}(A_{\alpha})&\text{if $x\in A_{\alpha}\in \mc A_{\eta}$ for some 
${\eta}<{\zeta}$},\\   
F(x)&\text{otherwise.}
 \end{array}
\right . 
\end{displaymath}
Then, by the inductive assumption $(\bullet_{|X_{\zeta}|})$ , there is a function
\begin{displaymath}
 c_{\zeta}: X_{\zeta}\to {\kappa}
\end{displaymath}
such that 
\begin{displaymath}
\text{ $c_{\zeta}(x)\notin F_{\zeta}(x)$ for $x\in X_{\zeta}$,
and $|I_{c_\zeta}(A_{\alpha})|={\kappa}$ 
for all $A_{\alpha}\in{\mc A}_{\zeta}$.
} 
\end{displaymath}
If $A_{\alpha}\in{\mc A}_{\zeta}$, then 
$|A_{\alpha}\cap X_{<{\zeta}}|<{\kappa}$ by (M2), and so
$|I_{c_\zeta}(A_{\alpha})|={\kappa}$  implies
\begin{displaymath}
|I_{c_{\le \zeta}}(A_{\alpha})|={\kappa}. 
\end{displaymath}

Since $\mc I$ is  a principal ideal, 
we can pick $G_{\zeta}(A_{\alpha})\in 
\mc I\cap \br I_{c_{\le \zeta}}(A_{\alpha});{\kappa};
$.

Finally take  $c=\cup \{c_{\zeta}:{\zeta}<\cf({\lambda})\}$.

Then  $c$ witnesses ($\bullet_{{\lambda}}$).
 ($\mc I$1) is clear from the construction.
As for  ($\mc I$2), assume that $A_{\alpha}\in {\mc A}_{\zeta}$.
Then $G_{\zeta}(A_{\alpha})\subset I_{c_{\le {\zeta}}}(A)$
by (\ref{eq:GsubsI}),
and if $x\in A_{\alpha}\setm X_{\le {\zeta}}$, then 
$c(x)\notin G_{\zeta}(A_{\alpha})$.
So  $G_{\zeta}(A_{\alpha})\subset I_{c}(A) $ as well.

So we completed the proof of Theorem \ref{tm:main}.
\end{proof}

\end{document}